\documentclass[a4paper,10pt]{amsart}
\usepackage{amsmath,amssymb,amsthm,latexsym}
\usepackage{verbatim}
\usepackage{mdwlist}
\usepackage[dvips]{graphicx}
\usepackage{psfrag}

\usepackage{changebar}

\newtheorem{thm}{Theorem}[section]
\newtheorem{lem}[thm]{Lemma}
\newtheorem{cor}[thm]{Corollary}
\newtheorem{prop}[thm]{Proposition}

\newtheorem{DEF}[thm]{Definition}
\theoremstyle{remark}
\newtheorem*{rem}{Remark}

\newcommand{\Hom}{\texttt{Hom}\,}
\newcommand{\Homc}{$\texttt{Hom}\,(C_m,C_n)$ }
\newcommand{\Homck}{$\texttt{Hom}\,(C_m,C_n)$}
\newcommand{\Homl}{$\texttt{Hom}\,(C_{2m},L_n)$ }
\newcommand{\Homlk}{$\texttt{Hom}\,(C_{2m},L_n)$}

\newcommand{\Hnulak}{$\texttt{Hom\,}_0(G,H)$}
\newcommand{\HomH}{$\texttt{Hom}\,(G,H)$ }
\newcommand{\HomHk}{$\texttt{Hom}\,(G,H)$}

\newcommand{\deset}{\vspace{10pt}}

\newcommand{\tn}{\textnormal}

\newcommand{\temesk}{$(a_1,a_2,\dots,a_m)$}
\newcommand{\temels}{$(a_1,a_2,\dots,a_{2m})$ }

\newcommand{\temes}{$(a_1,a_2,\dots,a_m)$ }
\newcommand{\teme}{$(i;i_1,\dots,i_r)$ }
\newcommand{\temelk}{$(i;i_1,\dots,i_m)$}
\newcommand{\temel}{$(i;i_1,\dots,i_m)$ }
\newcommand{\temek}{$(i;i_1,\dots,i_r)$}
\newcommand{\Cl}{\textnormal{Cl}}
\newcommand{\bdm}{\begin{displaymath}}
\newcommand{\edm}{\end{displaymath}}
\newcommand{\beq}{\begin{equation}}
\newcommand{\eeq}{\end{equation}}
\newcommand{\beqa}{\begin{eqnarray*}}
\newcommand{\eeqa}{\end{eqnarray*}}

\newcommand{\xt}{\widetilde{X}_i^m}
\newcommand{\pe}{\mathcal{P}}
\title[The homotopy type of complexes of graph homomorphisms]{The
homotopy type of complexes of graph homomorphisms between cycles}

\author{Sonja Lj. \v{C}uki\'{c}}
\address{Department of Mathematics, Royal Institute of Technology, Stockholm, Sweden\\ \indent current address: Department of Computer Science, Eidgen\"ossische Technische
Hoch\- schule, Z\"urich, Switzerland}
\email{sonja.cukic@inf.ethz.ch}
\author{Dmitry N. Kozlov}
\address{Department of Computer Science, Eidgen\"ossische Technische
Hochschule, Z\"urich, Switzerland}
\email{dkozlov@inf.ethz.ch}
\thanks {Research partially supported by Swedish Science Council}

\date{May, 2004}

\begin{document}
\begin{abstract}
In this paper we study the homotopy type of $\Hom(C_m,C_n)$, where
$C_k$ is the cyclic graph with $k$ vertices. We enumerate connected
components of $\Hom(C_m,C_n)$ and show that each such component is
either homeomorphic to a~point or homotopy equivalent to $S^1$.

Moreover, we prove that $\Hom(C_m,L_n)$ is either empty or is homotopy
equivalent to the union of two points, where $L_n$ is an $n$-string,
i.e., a tree with $n$ vertices and no branching points.
\end{abstract}

\maketitle

\section{Introduction}

To any two graphs $T$ and $G$ one can associate a~cell complex
$\Hom(T,G)$, see Definition~\ref{homdf}. The motivation for
considering $\Hom(T,G)$ came from the fact that it has good structural
properties, and that some special cases yield previously known
constructions. For example, $\Hom(K_2,G)$ is homotopy equivalent to
the neighborhood complex ${\mathcal N}(G)$, which plays the central
role in the Lov\'asz' proof of the Kneser Conjecture in 1978, see
\cite{Lov}.

On the other hand, since Babson \& Kozlov proved the Lov\'asz
Conjecture, \cite{BK3}, stating that for any graph $G$, and $r\geq 1$,
$k\geq -1$:
\begin{center}
{\it if }$\Hom(C_{2r+1},G)$ {\it is $k$-connected, then $\chi(G)\geq k+4$,}
\end{center}
it has become increasingly clear that the topology of $\Hom$-complexes
carries vital information pertaining to obstructions to the existence
of graph colorings. We refer the reader to the survey article \cite{IAS} for an introduction and further facts about \Hom-complexes.

Until now, the homotopy type of $\Hom(T,G)$ was computed only in very
few special cases. It was proved in \cite{CGH} that $\Hom(K_m,K_n)$ is
homotopy equivalent to a~wedge of $(n-m)$-dimensional spheres. It was
also shown that ``folding'' (i.e., removing a~vertex $v$ such that
there exists another vertex $u$ whose set of neighbors contains that
of $v$) the graph $T$ does not change the homotopy type of
$\Hom(T,G)$. This means, for example, that $\Hom(T,K_n)\simeq
S^{n-2}$, where $T$ is a~tree, since one can fold any tree to an~edge,
and since, as also shown in \cite{CGH}, $\Hom(K_2,K_n)$ is
homeomorphic to $S^{n-2}$. Beyond these, and a few other either degenerate
or small examples, nothing is known.

In this paper we study the homotopy type of the complex of graph
homomorphisms between two cycles, $\Hom(C_m,C_n)$, in particular
$\Hom(C_m,K_3)$, as well as between a cycle and a string. It is easy
to see that the connected components of $\Hom(C_m,C_n)$ can be indexed
by the signed number of times $C_m$ wraps around $C_n$ (with
an~additional parity condition if $m$ and $n$ are even). Also, if $n$
divides $m$, and $C_m$ wraps around $C_n$ $m/n$ times in either
direction, then, since there is no freedom to move, the corresponding
connected components are points. It was further noticed by the
authors, that in all up to now computed cases of $\Hom(C_m,C_n)$ it
turned out that all other connected components were homotopy
equivalent to~$S^1$.  The main result of this paper,
Theorem~\ref{mainth}, states that this is the case in general.

Our proof combines the methods of Discrete Morse theory with the
classical homotopy gluing construction, \cite[Section 4.G]{Hat}. In
order to be able to phrase our combinatorial argument concisely, we
develop a new encoding system for the cells of
$\Hom(C_m,C_n)$. Namely, we index the cells with collections of marked
points and pairs of points on circles of length~$m$, and translate the
boundary relation into this language.

The case $n=4$ is
a bit special and is dealt with separately, using the fact that the
folds are allowed in the second argument of $\Hom(-,-)$ as well, as
long as the removed vertex has an exact double in the set of the
remaining vertices, see Lemma~\ref{reduction}. The homotopy type of
the complex $\Hom(C_m,L_n)$ is computed by a~similar argument.

\section{Basic notations and definitions}

For any graph $G$, we denote the set of its vertices by $V(G)$, and 
the set of its edges by $E(G)$, where $E(G)\subseteq V(G) \times V(G)$. 
In this paper we will consider only undirected graphs, so $(x,y)\in E(G)$ 
implies that $(y,x)\in E(G)$. Also, our graphs are finite and may contain loops.\\
$\circ$ For a natural number $k$ we introduce the following notation 
$[k]=\{1,2,\dots,k\}$.\\
$\circ$ For a graph $G$ and $S\subseteq V(G)$ we denote by $G[S]$ the graph on the vertex set $S$ induced by $G$, that is $V(G[S])=S$, $E(G[S])=(S\times S)\cap E(G)$. We will denote the graph $G[V(G)\setminus S]$ by $G-S$.\\
$\circ$ Let $\texttt{N}(v)$ be the set of all neighbors of $v\in V(G)$, for a graph $G$, that is the set $\{w\in V(G)\: \vert\: (v,w)\in E(G)\}$.\\
$\circ$ For an integer $n\ge 2$, denote with $C_n$ and $L_n$, graphs such that $V(C_n)=V(L_n)=[n]$ and $E(C_n)=\{(x,x+1),(x+1,x)\,\vert\,x\in \mathbb{Z}_n\}$, $E(L_n)=\{(x,x+1),(x+1,x)\,\vert\,x\in [n-1]\}$.\\

\begin{DEF} For two graphs $G$ and $H$, a \textbf{graph homomorphism} 
from $G$ to $H$ is a map $\phi: V(G)\to V(H)$, such that if $x, y\in V
(G)$ are connected by an edge, then $\phi(x)$ and $\phi(y)$ are also
connected by an edge.\\ \indent We denote the set of all homomorphisms
from $G$ to $H$ by \textnormal{\Hnulak}.
\end{DEF}

\noindent The next definition is due to Lov\'{a}sz and was stated in
this form in \cite{TopO}.

\begin{DEF} \label{homdf}
\textnormal{\HomH} is a polyhedral complex whose cells are indexed by
all functions $\eta : V (G) \to 2^{V (H)}\setminus\{ \varnothing\}$,
such that if $(x, y)\in E(G)$, then for all $\tilde{x}\in \eta(x)$ and
$\tilde{y}\in \eta(y)$, $(\tilde{x},\tilde{y})\in E(H)$.

The closure of a cell $\eta$, $\Cl(\eta)$, consists of all cells
indexed by $\tilde{\eta}: V(G)\to 2^{V (H)}\setminus\{ \varnothing\}$
which satisfy the condition that $\tilde{\eta}(v)\subseteq {\eta}(v)$,
for all $v\in V(G)$.
\end{DEF}

It is easy to see that the set of vertices of \HomH is \Hnulak. Cells
of \HomH are direct products of simplices and dimension of a cell
$\eta$ is equal to $\sum_{v \in V(G)}\vert \eta (x) \vert - \vert
V(G)\vert $.\\

\begin{DEF} For an integer $i$, let $[i]_m$ be an integer such that $[i]_m\in[m]$ and $i\equiv [i]_m (\textnormal{mod } m)$.
\end{DEF}

In this paper we will deal mostly with \Homck. In this case each
vertex is denoted with $m$-tuple $(a_1,a_2,\dots,a_m)$, such that
$a_i\in [n]$ and $[ a_{[i+1]_m}-a_i]_n \in\{1,n-1\}$, for all
$i\in[m]$.  We also see that all cells of these complexes are cubes, since they are direct products of simplices and, clearly, the dimension of each simplex in this product is either 1 or 0.  \\

{\bf Some examples of \Homc complexes:}
\begin{itemize}
\item  \Homc is an empty set if $m$ is odd and $n$ is even: if $n$ is even, then there exists a map $\varphi:C_n\to K_2$, so if $\psi:C_m\to C_n$ exists then $\varphi\circ \psi:C_m\to K_2$ implying that $\chi(C_m)\leq 2$.  
\item Examples when $n=3$ and $m=2,4,5,6,7$ can be found in \cite{CGH}. The number of connected components of $\Hom (C_m,C_3)$ was computed there and it was proven that each connected component is either a point or homotopy equivalent to $S^1$. The following 3 cases were also observed in \cite{CGH}.
\item $\Hom(C_{2r+1},C_{2q+1})=\emptyset$ if and only if $r<q$.
 \item  $\Hom(C_{2r+1},C_{2r+1})$ is a disjoint union of 4r+2 points, for $r\geq 1$.
\item $\Hom(C_{2r+1},C_{2r-1})$ is a disjoint union of two cycles, length of each of them equal to $4r^2-1$.
\item $\Hom (C_4,C_{2r+1})$, for $r\geq 1$, is connected and it has $4r+2$ squares linked in the way depicted on Figure  \ref{homc4}.
\item $\Hom (C_4, C_{2r})$, for $r>2$, has two isomorphic connected components, each of them has $2r$ squares (see the Figure  \ref{homc4}).
\begin{figure}[ht] 
\begin{center} 
\psfrag{1}{$\scriptscriptstyle{(1,2,1,2)}$}
\psfrag{2}{$\scriptscriptstyle{(1,2,1,2r+1)}$}
\psfrag{3}{$\scriptscriptstyle{(1,2r+1,1,2r+1)}$}
\psfrag{4}{$\scriptscriptstyle{(1,2r+1,1,2)}$}
\psfrag{5}{$\scriptscriptstyle{(1,2r+1,2r,2r+1)}$}
\psfrag{6}{$\scriptscriptstyle{(1,2,3,2)}$}
\psfrag{Hom1}{$\Hom (C_4,C_{2r+1})$}
\psfrag{8}{$\scriptscriptstyle{(1,2,1,2)}$}
\psfrag{9}{$\scriptscriptstyle{(1,2,1,2r  )}$}
\psfrag{b}{$\scriptscriptstyle{(1,2r  ,1,2r  )}$}
\psfrag{a}{$\scriptscriptstyle{(1,2r  ,1,2)}$}
\psfrag{c}{$\scriptscriptstyle{(1,2r  ,2r-1,2r  )}$}
\psfrag{7}{$\scriptscriptstyle{(1,2,3,2)}$}
\psfrag{Hom2}{$\Hom (C_4,C_{2r})$}
\psfrag{e}{$\scriptscriptstyle{(2,1,2,1)}$}
\psfrag{f}{$\scriptscriptstyle{(2,1,2,3   )}$}
\psfrag{h}{$\scriptscriptstyle{(2,3   ,2,3   )}$}
\psfrag{g}{$\scriptscriptstyle{(2,3   ,2,1)}$}
\psfrag{i}{$\scriptscriptstyle{(2,3   ,4   ,3   )}$}
\psfrag{d}{$\scriptscriptstyle{(2,1,2r,1)}$}
\includegraphics[scale=0.6]{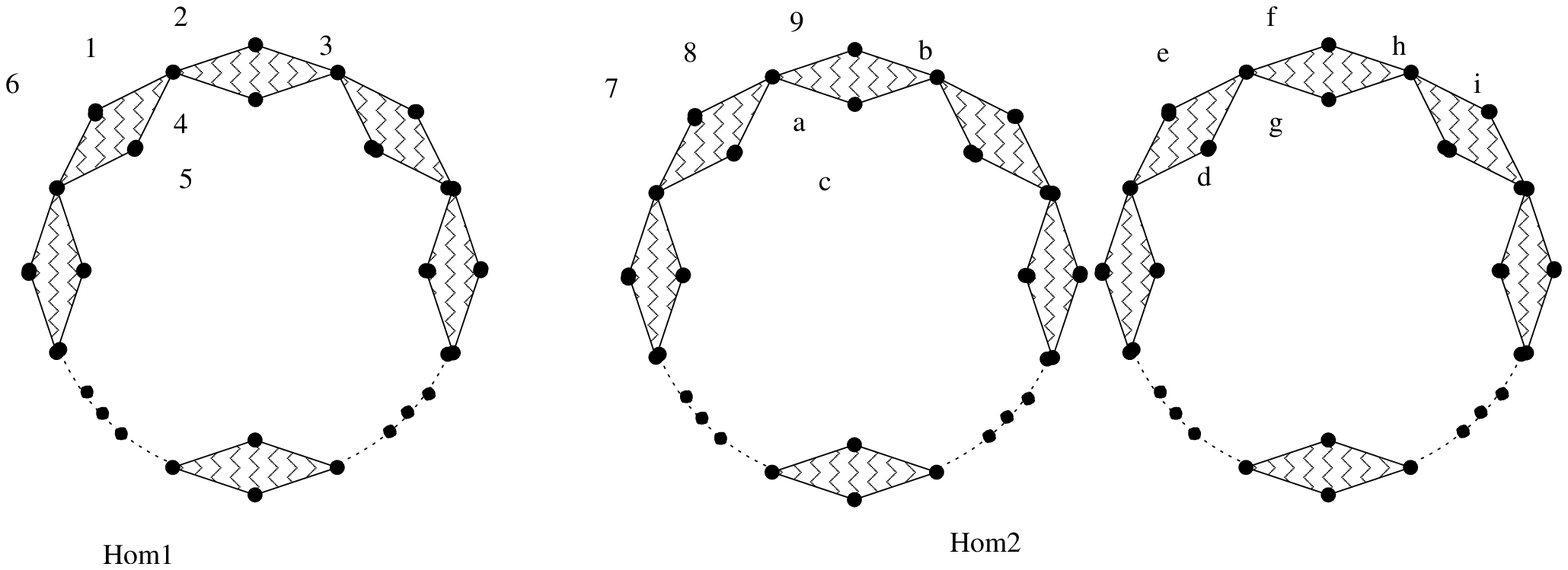}
\caption{} 
\label{homc4}
\end{center} 
\end{figure}
 \item $\Hom (C_9,C_3)$ consists of 6 isolated points  and two additional isomorphic connected parts, each of them has 90 solid cubes, 27 squares, 567 edges and 252 vertices. The local structure of one of those parts is shown on Figure \ref{c9c3}. The length of the cycle which is painted bold on the picture is 27.
\begin{figure}[ht] 
\begin{center} 
\includegraphics[scale=1]{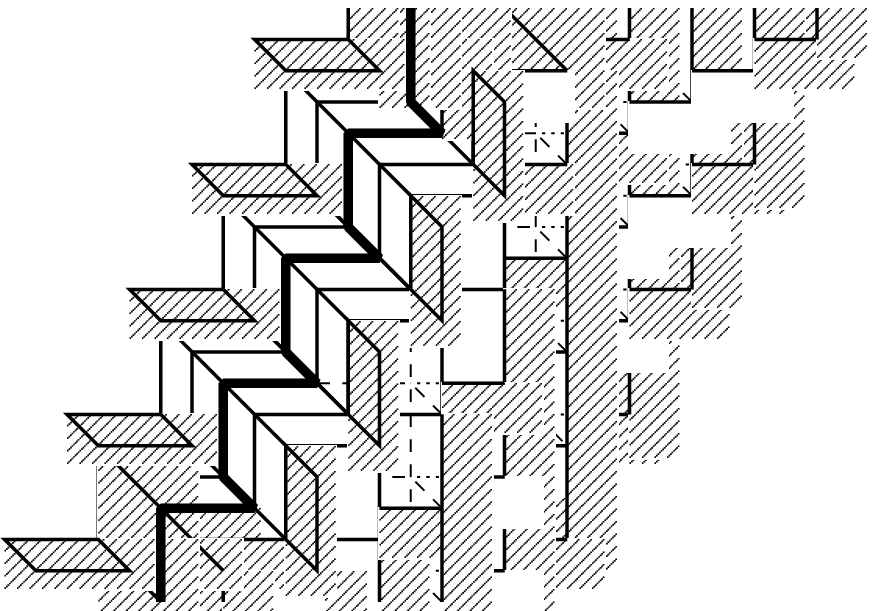}
\caption{} 
\label{c9c3}
\end{center} 
\end{figure}
\end{itemize}

\section{Complex {\textnormal{$\Hom(C_{2m},C_4)$}}}

We shall discuss this case separately because, unlike other cycles, $C_4$ has two vertices $u$ and $v$ such that $\texttt{N}(u)=\texttt{N}(v)$. In this section we will first prove a lemma for general $\Hom$ complexes, and then apply it to decide the homotopy type of  $\Hom(C_{2m},C_4)$.
\begin{lem}\label{reduction} Let $G$ and $H$ be graphs and $u,v\in V(H)$ such that \emph{$\texttt{N}(u)=\texttt{N}(v)$}. Also, let $i:H-v \hookrightarrow H$ be the inclusion and $\omega:H\to H-v$ the unique graph homomorphism which maps $v$ to $u$ and fixes other vertices. Then, these two maps induce homotopy equivalences \emph{$i_H:\Hom(G,H-v)\to \Hom(G,H)$} and \emph{$\omega_H:\Hom(G,H)\to \Hom(G,H-v)$}, respectively.
\begin{rem} A similar theorem about reduction of certain $\Hom$ complexes (Proposition 5.1) was proven in \cite{CGH}. This lemma was also proven independently in \cite{D}.
\end{rem}
\begin{proof} We will show that $\omega_H$ satisfies the conditions (A) and (B) of the Proposition 3.2 from \cite{CGH}.
Unfolding definitions, we see that for a cell of \HomHk, $\tau:V(G)\to 2^{V(H)}\setminus \{\emptyset\}$, we have 
\bdm
\omega_H(\tau)(x)=\left\{ \begin{array}{ll} \tau(x), & \textnormal{if }v\notin \tau (x); \\(\tau(x)\cup \{u\})\setminus \{v\}, & \textnormal{otherwise}. \end{array} \right.
\edm
Let $\eta$ be a cell of $\Hom (G,H-v)$, $\eta:V(G)\to 2^{V(H)\setminus \{v\}}\setminus \{\emptyset\}$. Then $\pe(\omega_H)^{-1}(\eta)$ is the set of all $\eta'$ such that, for all $x\in V(G)$,
\bdm
(*) \left\{ \begin{array}{ll}  \eta'(x)=\eta(x), & \textnormal{if }u\notin \eta(x); \\ \eta'(x)\cap \{u,v\}\neq \emptyset \textnormal{ and }\eta'(x)\setminus \{u,v\}=\eta(x)\setminus\{u\}, & \textnormal{otherwise.} \end{array} \right.
\edm 
It is easy to see that, because of the condition  $\texttt{N}(u)=\texttt{N}(v)$, all $\eta'$ satisfying $(*)$ belong to \HomHk. Take $\zeta\in \pe(\omega_H)^{-1}(\eta)$ such that 
\bdm
\zeta(x)=\left\{\begin{array}{ll} \eta(x), & u\notin \eta(x) \\ \eta(x)\cup \{v\}, & u\in \eta(x) \end{array}\right., \textnormal{ for all }x\in V(G).
\edm
Obviously, $\zeta$ is the maximal element of $\pe(\omega_H)^{-1}(\eta)$. It follows that $\Delta(\pe(\omega_H)^{-1}(\eta))$ is contractible and condition (A) is satisfied.\\

Take now any $\tau\in \pe(\omega_H)^{-1}(\Hom(G,H-v)_{\geq \eta})$. Then $\eta(x)\setminus\{u\}\subseteq \tau(x)\setminus \{u,v\}$ for all $x\in V(G)$ and, if $u\in \eta(x)$, then $\tau(x)\cap\{u,v\}\neq \emptyset$. The set $\pe(\omega_H)^{-1}(\eta)\cap \pe(\Hom(G,H))_{\leq \tau}$ consists of all cells $\eta'$ such that, for $x\in V(G)$, 
\bdm
\left\{ \begin{array}{ll}  \eta'(x)=\eta(x), & \textnormal{if }u\notin \eta(x); \\ \eta'(x)\cap \{u,v\}\neq \emptyset \textnormal{ and } \eta(x)\setminus \{u\} \subseteq \eta'(x)\subseteq (\eta(x)\cup\{v\})\cap \tau(x), & \textnormal{otherwise,} \end{array} \right.
\edm 
and hence has the maximal element $\xi$, where  $\xi(x)=\eta(x)$ for $x\in V(G)$ such that $u\notin \eta(x)$ and $\xi(x)= \tau(x)\cap (\eta(x)\cup \{v\})$ otherwise. \\
Since it satisfies conditions (A) and (B), we conclude that $\textnormal{Bd}(\omega_H)$ and hence also $\omega_H$ are homotopy equivalences. \\

It is left to prove that $i_H$ is also a homotopy equivalence. It is clear that $\omega_H\circ i_H=id_{\Hom(G,H-v)}$. Let $\vartheta$ be the homotopy inverse of $\omega_H$. Then we have $i_H\circ \omega_H \simeq \vartheta \circ \omega_H \circ i_H \circ \omega_H \simeq  \vartheta \circ \omega_H \simeq id_{\Hom (G,H)}$. 
\end{proof}
\end{lem}

Now we have everything we need to prove the following theorem.
\begin{thm} \label{l2l3} The complex \emph{$\Hom (C_{2m},C_4)$}, for $m\geq 1$, is homotopy equivalent to a complex consisting of two points.
\begin{proof} We use Lemma \ref{reduction} and obtain 
\bdm
\Hom (C_{2m},C_4)\simeq \Hom (C_{2m},L_3) \simeq \Hom (C_{2m},L_2)
\edm
It is trivial to see that $\Hom (C_{2m},L_2)$ has two vertices, namely $(1,2,1,2,\dots,1,2)$ and $(2,1,2,1,\dots,2,1)$, and no other cells.
\end{proof}
\end{thm}
\section{Discrete Morse theory}
In this section we will introduce the notations and state the reformulation of Forman's result from Discrete Morse Theory given in \cite{CGH}. For more general results about this topic, see \cite{For}.
\begin{DEF} A \emph{partial matching}  on a poset $P$ with covering relation $\succ$ is a set $S\subseteq P$ together with an injective map $\mu: S\to P\setminus S$ such that $\mu(x)\succ x$, for all $x\in S$.
The elements from $P\setminus (S\cup \mu (S))$ are called \emph{critical}.
\end{DEF}

\begin{DEF} A matching is called \emph{acyclic} if there does not exist a sequence $x_0,x_1,\dots,x_t=x_0\in S$ such that $x_0\neq x_1$ and $\mu(x_i)\succ x_{i+1}$, for $i\in [t-1]$.
\end{DEF}
\noindent For a regular CW complex $X$ let $\mathcal{P}(X)$ be its face poset with covering relation $\succ$.
\begin{prop}\label{morse} Let $X$ be a regular CW complex, and let $X'$ be a subcomplex of $X$. Then the following are equivalent: 
\begin{enumerate}
\item there is a sequence of collapses leading from $X$ to $X'$;
\item there is an acyclic partial matching $\mu$ on $\mathcal{P}(X)$ with the set of critical cells being $\mathcal{P}(X')$. 
\end{enumerate}
\end{prop}
For proof see \cite[Proposition 5.4]{RatHom}.

\section{Another notation for the cells of {\textnormal{$\Hom (C_m,C_n)$}}}

\begin{rem} From now on, unless otherwise stated, we will work only with $\Hom (C_m,C_n)$ where $n\neq 4$.
\end{rem}
\begin{DEF} We say that $i\in [m]$ is \emph{a returning point} of a vertex \temes of \emph{\Homc} if  $[a_i-a_{[i+1]_m}]_n =1$.
\end{DEF}  
We see that each vertex \temes$\in \Hom _0 (C_m,C_n)$ uniquely determines an $(r+1)$-tuple $(i;i_1,\dots,i_r)$, where $m=nk+2r$, $i=a_1$ and $i_1,\dots,i_r$ are all its returning points with condition that $i_1<i_2<\dots<i_r$. Conversely, assume that we have a $(\rho+1)$-tuple $(j;j_1,\dots,j_\rho)$, where $0\leq \rho\leq m$, $j\in [n]$, $1\leq j_1<j_2<\dots<j_\rho \leq m$ and $m=kn+2\rho$, for some integer $k$. Then \temesk, defined by $a_i=[j+i-1-2\rho_i]_n$, where $\rho_i=\big\vert\{q\vert j_q<i\}\big\vert$, is a vertex of \Homc with returning points $j_1,\dots,j_\rho$. Indeed, $[a_m-a_1]_n=[j+m-2\rho\pm 1-j]_n=[kn\pm 1]_n\in \{1,n-1\}$, where $[a_m-a_1]_n=1$ if and only if $m=j_\rho$; and for all $i\in [m-1]$, $[a_{i}-a_{i+1}]_n=1$ if and only if $i=j_q$ for some $q\in [\rho]$, otherwise $[a_{i}-a_{i+1}]_n=n-1$. Hence, we have proven the following lemma:
\begin{lem} Let $S$ be a set containing all $(r+1)$-tuples \temek, such that $0\leq r\leq m$ and, for some integer $k$, $m=nk+2r$, $i\in [n]$ and $1\leq i_1<\dots<i_r\leq m$. Then there is a bijection $\Xi$ between \emph{$\Hom_0(C_m,C_n)$} and $S$ given by 
\bdm
\Xi((a_1,a_2,\dots,a_m))=(i;i_1,\dots,i_r), 
\edm
where $i=a_1$ and $i_1,\dots,i_r$ are all returning points of \temesk.
\end{lem}  

\begin{figure}[ht] 
\begin{center} 
\psfrag{1}{1}
\psfrag{2}{2}
\psfrag{5}{3}
\psfrag{Hom1}{$m=6,n=9$}
\psfrag{v1}{$(1,9,8,9,1,9)$}
\psfrag{Hom2}{$m=9,n=5$}
\psfrag{Hom3}{$m=8,n=8$}
\psfrag{v2}{$(3,4,5,1,5,1,2,1,2)$}
\psfrag{v3}{$(2,1,8,7,8,1,2,1)$}
\includegraphics[scale=1]{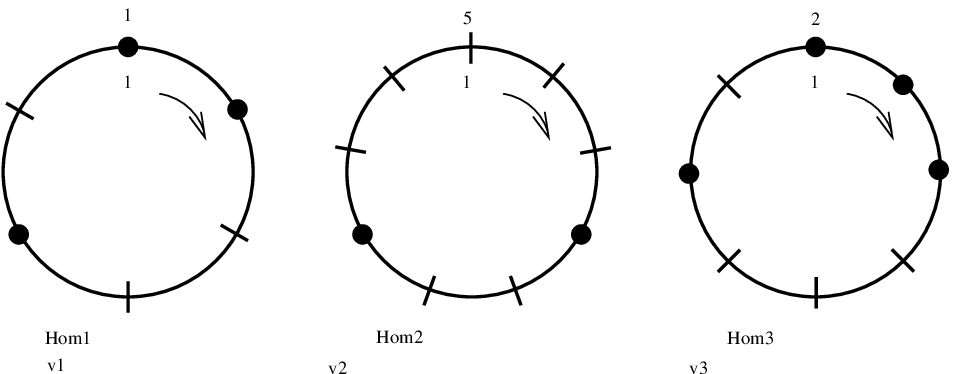}
\caption{} 
\label{uvod_cikl}
\end{center} 
\end{figure}

Sometimes we will represent a vertex \teme of \Homc by a picture of $C_m$ with emphasized returning points and number $i$. Some examples of such representation are shown on Figure \ref{uvod_cikl}. Vertices of $C_m$ are always ordered like this: if we start from vertex labeled with $1$ and go in clockwise direction, we will get an increasing sequence of numbers from $1$ to $m$.
\small

\normalsize
\begin{lem}\label{conn_comp} If two vertices are in the same connected component of \emph{\Homck}, then their number of returning points is the same.
\begin{proof} It is enough to prove that for an arbitrary edge, its endpoints have the same number of returning points.\\ Let $(a_1,\dots,a_i,\{[a_i-1]_n,[a_i+1]_n\},a_i,a_{i+3},\dots,a_{m})$ be an edge, and let $x$ and $y$ be its endpoints, $x=(a_1,\dots,a_i,[a_i-1]_n,a_i,a_{i+3},\dots,a_{m})$ and $y=(a_1,\dots,a_i,[a_i+1]_n,a_i,a_{i+3},\dots,a_{m})$. Now it is trivial to see that, for $j\in [m]\setminus\{i,i+1\}$, $j$ is a returning point for $x$ if and only if it is a returning point for $y$. Also, we see that $i$ is a returning point for $x$, while $i+1$ is not and, similarly, $i+1$ is a returning point for $y$ and $i$ is not. Hence, $x$ and $y$ have the same number of returning points. 
\end{proof}
\end{lem} 

 It follows from definitions that each cell $\eta\in$\Homc has the property that for any $x\in [m]$, $\vert \eta(x) \vert \in \{1,2\}$ and $\eta(x)$, $\eta([x+1]_m)$ cannot both have cardinality two. Also, if $\vert \eta(x) \vert =2$, then for some $i\in[n]$, $\eta([x-1]_m)=\eta([x+1]_m)=\{i\}$ and $\eta(x)=\{[i-1]_n,[i+1]_n\}$. \\ 

Now, let $\eta\in$\Homc be a cell. Since $\Cl(\eta)$ is connected, by Lemma \ref{conn_comp}, all its vertices have the same number of returning points. Denote that number with $r$. Then we can denote the cell $\eta$ with $(r+1)$-tuple of symbols $(s;s_{i_1},\dots,s_{i_r})$, where for all $j\in [r]$, $i_j\in [m]$, $1\leq i_1<\dots<{i_r}\leq m$ and
\bdm
s_{i_k}=\left\{\begin{array}{ll} i_k^+, & \textnormal{if }\vert \eta([i_k+1]_m)\vert =2; \\
i_k, & \textnormal{if } \eta(i_k)=\{j\} \textnormal{ and } \eta([i_k+1]_m)=\{[j-1]_n\} \textnormal{ for some }j\in [n]. \end{array} \right.
\edm 
Also, if $\eta(1)=\{i,[i+2]_n\}$ or if $\eta(1)=\{i\}$, for some $i\in [n]$, then $s=i$. Note that, if $s_{i_k}=i_k^+$, then  $[i_k+1]_m\neq i_{[k+1]_r}$\\
It is clear that each cell uniquely determines such an $(r+1)$-tuple of symbols and that dimension of a cell is equal to number of $s_{i_k}$ such that $s_{i_k}=i_k^+$.\\

Conversely, if $(j;s_{j_1},\dots,s_{j_r})$ is a $(r+1)$-tuple of symbols and if following conditions are satisfied:
\begin{enumerate}
\item $0\leq r\leq m$, $j\in [n]$, $m=nk+2r$ and $1\leq j_1 <j_2<\dots<j_r\leq m$,
\item For all $k\in[r]$, $s_{j_k}\in\{j_k,j_k^+\}$,
\item If $s_{j_k}=j_k^+$ then $[j_k+1]_m\neq j_{[k+1]_r}$,
\end{enumerate}
then it is not hard to check that $(j;s_{j_1},\dots,s_{j_r})$ corresponds exactly to one cell from \Homck, namely to $(A_1,\dots, A_m)$, where for $a_k=j+k-1-2\big\vert\{q\vert j_q<k\}\big\vert$, 
\bdm 
A_k=\left\{ \begin{array}{ll} \{ [a_k]_n\}, & \textnormal{if }[k-1]_m\notin \{i_l\vert s_{i_l}=i_l^+\}; \\  \{ [a_k]_n,[a_k+2]_n\}, & \textnormal{otherwise}. \end{array} \right.
\edm
   
\begin{figure}[ht] 
\begin{center} 
\psfrag{1}{1}
\psfrag{2}{2}
\psfrag{5}{5}
\psfrag{Hom1}{$m=6,n=9$}
\psfrag{v1}{$(1,9,\{1,8\},9,1,9)$}
\psfrag{u1}{$(1;1,2^+,5)$}
\psfrag{Hom2}{$m=9,n=5$}
\psfrag{Hom3}{$m=8,n=8$}
\psfrag{v2}{$(5,1,2,3,\{2,4\},3,4,\{3,5\},4)$}
\psfrag{u2}{$(5;4^+,7^+)$}
\psfrag{v3}{$(\{2,4\},3,2,\{1,3\},2,3,4,3)$}
\psfrag{u3}{$(2;2,3^+,7,8^+)$}
\includegraphics[scale=1]{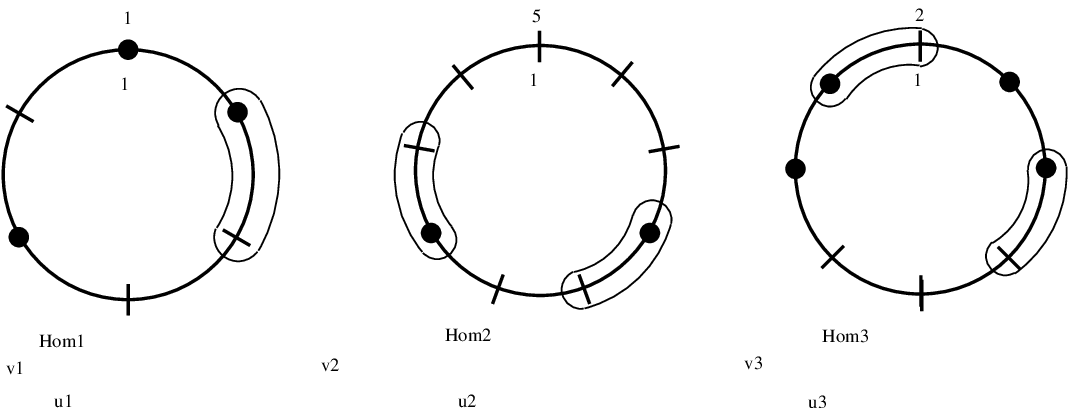}
\caption{} 
\label{cells_uvod}
\end{center} 
\end{figure}

 We will also represent those $(r+1)$-tuples (and corresponding cells) with pictures (see Figure \ref{cells_uvod}).
\small

\normalsize
\begin{rem} Assume $\eta=(s;s_{i_1},\dots,s_{i_r})\in$\Homck, and $s_{i_k}=i_k^+$, for $i_k\neq m$. "Unplusing" $\eta$ in the $k$-th position yields two  $(\dim(\eta)-1)$-dimensional cells in  $\Cl(\eta)$: $(s;s_{i_1},\dots,s_{i_{k-1}},i_k,s_{i_{k+1}},\dots,s_{i_r})$ and $(s;s_{i_1},\dots,s_{i_{k-1}},i_k+1,s_{i_{k+1}},\dots,s_{i_r})$.
\end{rem}
\begin{lem}\label{cls}  Let $\eta=(i;s_{i_1},\dots,s_{i_r})$ be a cell from \emph{\Homc} of dimension $d\geq 1$ and let $D=\{k\vert k\in[r]\textnormal{ and }s_{i_k}=i_k^+\}$. For $  k\in D\setminus \{m\}$, let 
\beqa
\eta_k &=& (i;s_{i_1},\dots,s_{i_{k-1}},i_k,s_{i_{k+1}},\dots,s_{i_r})\\
\eta_k'&=& (i;s_{i_1},\dots,s_{i_{k-1}},i_k+1,s_{i_{k+1}},\dots,s_{i_r}).
\eeqa
Then the set of all cells of \emph{\Homc} which are contained in \emph{$\Cl(\eta)$} and have dimension $d-1$ is equal to:\\
\tn{(1)} $\bigcup_{k\in D} \{\eta_k,\eta_k'\}$, if $m\notin D$,\\
\tn{(2)} $\bigcup_{\begin{subarray} {c} k\in D \\ k\neq m \end{subarray}} \{\eta_k,\eta_k'\}\cup \{(i;s_{i_1},\dots,s_{i_{r-1}},m),([i+2]_n;1,s_{i_1},\dots,s_{i_{r-1}})\}$, if $m\in D$.
\begin{proof} If we write the cell $\eta$ using old notation, it is clear that $\Cl(\eta)$ contains exactly $2d$ cells of dimension $d-1$. By the previous remark, we know that $S=\bigcup_{k\in D\setminus\{m\}} \{\eta_k,\eta_k'\}$ is a set consisting of different $(d-1)$-dimensional cells of $\Cl(\eta)$. \\
\indent In the case (1), $\vert S \vert=2d$ and $S$ is exactly the set of all $(d-1)$-dimensional cells in $\Cl(\eta)$.  \\
\indent Let us now deal with the case (2). Then $\eta=(\{i,[i+2]_n\},[i+1]_n,\dots,[i+1]_n)$. Two of $(d-1)$-dimensional cells in $\Cl(\eta)$ are $\eta'= (i,[i+1]_n,\dots,[i+1]_n)$ and $\eta''=([i+2]_n,[i+1]_n,\dots,[i+1]_n)$ (we have changed only the first component). It is easy to see that, in new notations, $\eta'=(i;s_{i_1},\dots,s_{i_{r-1}},m)$ and $\eta''=([i+2]_n;1,s_{i_1},\dots,s_{i_{r-1}})\}$. The claim follows since $\vert S\cup \{\eta',\eta''\} \vert =2d$ and all cells are different.
\end{proof}
\end{lem}
\begin{rem} From our definition of $\Cl(\eta)$ we see that $\Cl(\eta)\setminus\{\eta\}=\cup_{\eta'\in S} \Cl (\eta')$, where $S$ is the set of  all cells of dimension $\dim(\eta)-1$ contained in $\Cl(\eta)$.
\end{rem}

\begin{figure}[ht] 
\begin{center} 
\small 
\psfrag{dva}{2}
\psfrag{jed}{1}
\psfrag{(2;2,3^+,6^+,8^+)}{$(2;2,3^+,6^+,8^+)$}
\psfrag{(2;2,3,6^+,8^+)}{$\scriptstyle{(2;2,3,6^+,8^+)}$}
\psfrag{(2;2,3^+,6,8^+)}{$\scriptstyle{(2;2,3^+,6,8^+)}$}
\psfrag{(2;2,3^+,6^+,8)}{$\scriptstyle{(2;2,3^+,6^+,8)}$}
\psfrag{(2;2,4,6^+,8^+)}{$\scriptstyle{(2;2,4,6^+,8^+)}$}
\psfrag{(2;2,3^+,7,8^+)}{$\scriptstyle{(2;2,3^+,7,8^+)}$}
\psfrag{(4;1,2,3^+,6^+)}{$\scriptstyle{(4;1,2,3^+,6^+)}$}
\psfrag{1}{$\scriptstyle{1} $}
\psfrag{2}{$\scriptstyle{2} $}
\psfrag{4}{$\scriptstyle{4} $}
\includegraphics[scale=1]{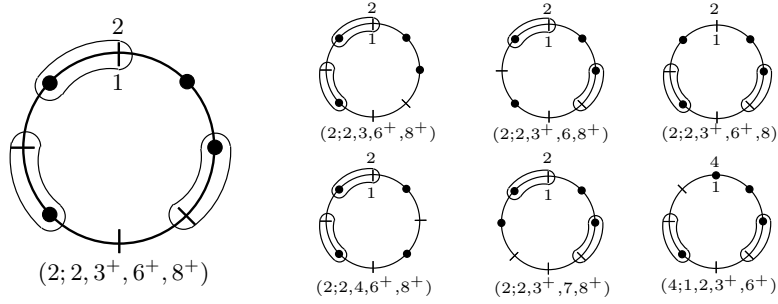}
\caption{A cell from $\Hom (C_8,C_8)$ of dimension 3 and all 2-dimensional cells contained in its closure.} 
\label{subcells}
\end{center} 
\end{figure} 
\begin{lem}\label{conn_comp2} Suppose $[r]_m\neq m$. The two vertices \teme and $(j;j_1,\dots,j_r)$  of \emph{\Homc} are in the same connected component if and only if $[i+2l]_n=j$, for some integer $l$.
\begin{proof}
Suppose first that $[i+2l]_n=j$. \\
By the previous remark, if $k\neq r,i_{k+1}>i_k+1$ or $k=r,i_k\neq m$, then vertices \teme and $(i;i_1,\dots,i_{k-1},i_k+1,i_{k+1},\dots, i_r)$ are in the same connected component (both of them are elements of $\Cl (i;i_1,\dots,i_{k-1},i_k^+,i_{k+1},\dots, i_r)$). Let us introduce an equivalence relation $\sim$ on the set $\texttt{Hom}\,_0(C_m,C_n)$ as follows: $x\sim y$ if and only if $x$ and $y$ lie in the same connected component. Using these two things, it is easy to see that:
\bdm
\begin{split}
(i;i_1,\dots,i_{r-1},i_r) &\sim (i;i_1,\dots,i_{r-1},m)\sim (i;i_1,\dots,m-1,m)\sim \cdots \sim \\ &\sim(i;m-r+1,\dots,m-1,m) 
\end{split}
\edm
Now we have (see the proof of Lemma \ref{cls} and use the fact that $[r]_m\neq m$):
\bdm
\begin{split}
 (i;i_1,\dots,i_r) &\sim (i;m-r+1,\dots,m-1,m) \sim ([i+2]_n;1,m-r+1,\dots,m-1)\\ &\sim ([i+2]_n;m-r+1,\dots,m-1,m) \sim \dots \sim \\ &\sim ([i+2l]_n;m-r+1,\dots,m-1,m) \sim  (j;j_1,\dots,j_r).
\end{split}
\edm
Conversely, suppose that for all integers $l$, $[i+2l]_n\neq j$. This can happen only if $n$ is even number and $[i]_2\neq[j]_2$. Since the parity of the first coordinate is constant on edges, we see that \teme and $(j;j_1,\dots,j_r)$ cannot be in the same connected component.
\end{proof}
\end{lem}
\begin{rem} If \temes is a vertex with $r$ returning points, where $[r]_m=m$, then it is easy to see that for all $i\in m$, $a_{[i+2]_m}\neq a_i$. Hence, there does not exist an edge with this vertex as an endpoint. 
\end{rem}

\normalsize
\section{The homotopy type of {\textnormal{\Homc}}}

\begin{thm} \label{mainth} Assume \emph{\Homc}$\neq\emptyset$, and let $X$ be some connected component of \emph{\Homc}. Then $X$ is either a point or is homotopy equivalent to a circle.
\begin{proof}
By Theorem \ref{l2l3} we know that the statement is true for $n=4$, so assume $n\neq 4$. \\ Lemma \ref{conn_comp}  implies that all vertices in $X$ have the same number of returning points. Denote that number with $r$. If $r=0$ or $r=m$, then, clearly, $X$ is a point. Let us deal with the case when $0<r<m$.\\
For all $i\in [n]$, let $X_i$ be the subcomplex of $X$ consisting of closures of all cells $\eta$ such that $i\in \eta(1)$ and let $\xt$ be the induced subcomplex of $X$ on the vertices $\eta$ such that $\eta(1)=\{i\}$. It is obvious that $\xt \subseteq X_i$ and that $X=\bigcup_{i=1}^n X_i$. Notice that, in the case when $n$ is even, it is a corollary of Lemma \ref{conn_comp2} that either $X_1=X_3=\cdots=X_{n-1}=\emptyset$ or $X_2=X_4=\cdots=X_n=\emptyset$. \\

\noindent {\it Claim 1.} $\xt$ is a strong deformation retract of $X_i$, for all $i$ such that $X_i\neq \emptyset$.
\begin{proof}[Proof of Claim 1]
Let us define a partial matching on $\mathcal{P}(X_i)$ on the following way: for $\eta\in \mathcal{P}(X_i)$ such that $i \notin \eta(1)$, let $\mu(\eta):=\widetilde{\eta}$ where
\bdm
\widetilde{\eta}(j)=\left\{ \begin{array}{ll} \eta(1)\cup \{i\}, & \textnormal{for }j=1; \\ \eta(j), & \textnormal{for }j=2,3,\dots,m.  \end{array} \right. 
\edm 
Obviously this is an acyclic matching and $\eta$ is a critical cell if and only if $\eta(1)=\{i\}$. Hence all critical cells form the subcomplex $\xt$. By Proposition \ref{morse}, there exists a sequence of collapses from $X_i$ to $\xt$, and since a collapse is a strong deformation retract, we see that the claim 1 is true and $\xt$ and $X_i$ have the same homotopy type.
\end{proof}
\emph{Remark.} It is important to see that there \emph{does not} exist an edge between vertices $v=(i;1,i_2,\dots,i_r)$ and $w=(i;i_2,\dots,i_r,m)$ since, in the old notation, $v=(i,i-1,\dots,i-1)$ and $w=(i,i+1,\dots,i+1)$ and two vertices of the same edge can be different only on one coordinate. Because of that, when we want to give a picture for a better explanation, instead of $C_m$ we will draw $L_m$. The reason for this notation is that there is no $m^+$ in notations of cells of $\xt$. For example, the cell $(4;1,2,3^+,6)$ from $\Hom (C_8,C_8)$ we will represent as in Figure \ref{prim}.
\begin{figure}[ht] 
\begin{center}  
\psfrag{1}{$\scriptstyle{1} $}
\psfrag{4}{$\scriptstyle{4} $}
\includegraphics[scale=1]{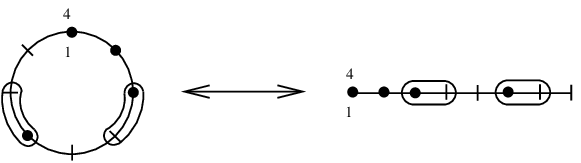}
\caption{} 
\label{prim}
\end{center} 
\end{figure} 

\noindent {\it Claim 2.} If $\xt\neq \emptyset$, then it is contractible.
\begin{proof}[Proof of Claim 2]
Before we define a matching on $\mathcal{P}=\mathcal{P}(\xt)$, we will need some additional notations and definitions. For a cell $\eta=(i;s_{i_1},\dots,s_{i_r})$ let $F(\eta)= \{i_j\vert j\in[r],i_j\neq m,s_{i_j}=i_j\textnormal{ and }i_j+1\notin \{i_1,\dots,i_r\}\}$ and $F^+(\eta)= \{i_j\vert j\in[r],s_{i_j}=i_j^+\}$. Also, we define two maps $R,R^+:\xt \to \mathbb{N}\cup\{\infty\}$ in the following way:\\ For a cell $\eta=(i;s_{i_1},\dots,s_{i_r})$ let
\bdm
R(\eta)=\left\{\begin{array}{ll}  \infty, & F(\eta)=\emptyset \\
\min  F(\eta), & \textnormal{otherwise} \end{array} \right.,
\textnormal{ and }
R^+(\eta)=\left\{\begin{array}{ll}  \infty, & F^+(\eta)=\emptyset \\
\min  F^+(\eta), & \textnormal{otherwise} \end{array} \right..
\edm  
Now, let $S=\{\eta\in \mathcal{P} \vert R(\eta)\lneq R^+(\eta)\}$. In particular, if $\eta\in S$, then $R(\eta)\neq \infty$.  For $\eta=(i;s_{i_1},\dots,s_{i_r})\in S$, let $\nu(\eta)=(i;p_{i_1},\dots,p_{i_r})$, where
\bdm
p_{i_j}=\left\{\begin{array}{ll} s_{i_j}, & i_j\neq R(\eta); \\ i_j^+, & \textnormal{otherwise.} \end{array} \right.
\edm

\small
\begin{figure}[ht] 
\begin{center}  
\psfrag{1}{$\scriptstyle{1} $}
\psfrag{i}{$\scriptstyle{i} $}
\psfrag{2}{$\nu(\eta_b)$}
\psfrag{2'}{$\eta_b$}
\psfrag{3}{$\nu(\eta_c)$}
\psfrag{3'}{$\eta_c$}
\psfrag{4}{$\nu(\eta_a)$}
\psfrag{4'}{$\eta_a$}
\includegraphics[scale=1]{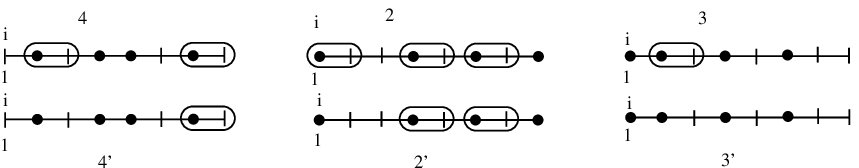}
\caption{} 
\end{center} 
\end{figure} 

\normalsize
It is clear that $\nu$ is injective, $\nu(\eta)\succ \eta$  and $R^+(\nu(\eta))<R(\nu(\eta))$, hence $\nu(\eta)\notin S$. We conclude that $(S,\nu)$ is a partial matching on $\mathcal{P}$. We will now prove that $\sigma=(i;m-r+1,m-r+2,\dots,m)$ is the only critical cell of our matching. Let $(i;s_{i_1},\dots,s_{i_r})=\xi \in \mathcal{P}\setminus S$. Then we have two cases: 
\begin{enumerate}
\item $\infty=R(\xi)=R^+(\xi)$, then $\xi=\sigma$;
\item $R(\xi)>R^+(\xi)$. In this case, let $\xi'=(i;p_{i_1},\dots,p_{i_r})$, where 
\bdm
p_{i_j}=\left\{\begin{array}{ll} s_{i_j}, & i_j\neq R^+(\eta); \\ i_j, & \textnormal{otherwise.} \end{array} \right.
\edm
Obviously, $\xi'\in S$, $\nu(\xi')=\xi$, hence $\xi$ is not a critical cell.
\end{enumerate}
 We conclude that $\mathcal{P}\setminus (S\cup \nu(S))=\{\sigma\}$. \\

What is left to prove is that this matching is acyclic. For any cell  $\eta\in \xt$, $\eta=(i;s_{i_1},\dots,s_{i_r})$, let $\Sigma(\eta)=\sum_{j=1}^r i_j$. Notice that $\Sigma(\eta)=\Sigma(\nu(\eta))$, for $\eta \in S$. \\
If $\xi=(i;l,l+1,\dots,l+k,s_{i_{k+2}},\dots,s_{i_r})\in S$, where $R(\xi)=l+k$, and if $\xi'\neq \xi$ is a cell such that $\nu(\xi)\succ \xi'$, it is not hard to see that, if $\xi' \in S$, then we must have $\xi'=(i;l,l+1,\dots,l+k-1,l+k+1,s_{i_{k+2}},\dots,s_{i_r})$. Then $\Sigma(\xi')=\Sigma(\xi)+1$ (here we have also used the previous remark). \\ Hence, if $\eta_1,\dots,\eta_t\in S$ such that $\eta_1\neq\eta_2$ and $\nu(\eta_i)\succ \eta_{i+1}$, then $\Sigma(\eta_t)>\Sigma(\eta_1)$ and it is not possible that $\nu(\eta_t)\succ \eta_1$, since in that case it would have to be $\Sigma(\eta_1)=\Sigma(\eta_t)+1$. 

\small{
\begin{figure}[ht] 
\begin{center}  
\psfrag{1}{$\eta_1$}
\psfrag{2}{$\eta_2$}
\psfrag{3}{$\eta_3$}
\psfrag{4}{$\eta_4$}
\psfrag{5}{$\eta_5$}
\psfrag{6}{$\eta_6$}
\includegraphics[scale=1]{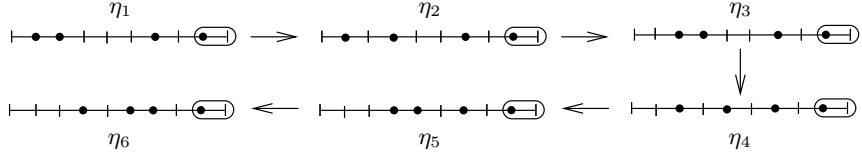}
 \caption{Example of a sequence $\eta_1,\dots,\eta_6\in S$ such that $\nu(\eta'_i)\succ \eta_{i+1}$, showing how returning points are ``moving'' to the right}
\end{center} 
\end{figure}
} 
\normalsize

By Proposition \ref{morse} there exists a sequence of elementary collapses leading from $\xt$ to $\{\sigma\}$ and hence $\xt$ is contractible.
\end{proof}

 We have proven that all non-empty subcomplexes $X_i$ are contractible. Let us now determine the structure of their intersections. \\
Let $\eta=(A_1,A_2,\dots,A_m)$ be a cell such that $\eta \in X_i\cap X_j$, where $i\neq j$. Then, because of the definition of $X_i$ and $X_j$, we have that $(A_1\cup \{i\},A_2,\dots,A_m)\in X_i$ and $(A_1\cup \{j\},A_2,\dots,A_m)\in X_j$. But then, for all $k\in A_2\cup A_m$, it must be $[k-i]_n,[k-j]_n \in \{1,n-1\}$. This is possible only when $j=[i\pm 2]_n$ (since $n\neq 4$). Without any loss of generality we can assume that $j=[i+2]_n$ and then $A_2=A_m=\{[i+1]_n\}$.  We conclude that: 
\bdm
X_i\cap X_j=\left\{ \begin{array}{ll} I \times \widetilde{X}_{[i+1]_n}^{m-2}, & \textnormal{for }j=[i+2]_n\textnormal{ and }X_i\neq \emptyset; \\ I \times \widetilde{X}_{[i-1]_n}^{m-2}, & \textnormal{for }j=[i-2]_n\textnormal{ and }X_i\neq \emptyset; \\ \emptyset, & \textnormal{otherwise};  \end{array} \right. 
\edm 
where $\widetilde{X}_{[i\pm1]_n}^{m-2}$ are both subcomplexes of $\Hom(C_{m-2},C_n)$, and $I$ is the unit interval. Notice that  each vertex of $\widetilde{X}_{[i+1]_n}^{m-2}$ has $r-1\geq0$ returning points. We know that $\widetilde{X}_{[i\pm1]_n}^{m-2}$ are contractible, and hence $X_i\cap X_j$ is also contractible. We also see that $X_{k_1}\cap\dots\cap X_{k_t}=\emptyset$ for $t\geq 3$ (using the similar argument as in case $X_i\cap X_j$).

The family of subcomplexes $\{X_i\}_{i=1}^n$ satisfies the conditions of
\cite[Corollary 4G.3, Exercise 4G.4]{Hat}, and
hence we have $X\simeq \mathcal{N}(X_i)$, where $\mathcal{N}(X_i)$ is
the nerve of $\{X_i\}_{i=1}^n$. On the other hand,
$\mathcal{N}(X_i)=\{i,\{i,[i+2]_n\}\,\vert i\in [n] \textnormal{ and
}X_i\neq \emptyset \}\simeq S^1$ and therefore $X\simeq S^1$.
\end{proof}
\end{thm}

Now it is time to summarize our results for complexes \Homc. We will have several cases depending on the parity of both $m$ and $n$. First of all, let us specify how we will index different connected components. By Lemma \ref{conn_comp2} we have two essentially different cases:
\begin{itemize}
\item $n$ is even: We will denote with $\Delta_r^i$ connected component in which all vertices have $r\neq 0,\;m$ returning points and for some of its vertices $\eta$, $\eta(1)=\{i\}$, where $i\in\{1,2\}$.
\item $n$ is odd: In this case the only thing which determines a connected component is number of returning points $r\neq 0,\; m$ of any vertex in that component. Hence, we will denote it with $\Delta_r$.
\end{itemize}
In Table 9 we will denote with $[m]$ ($[m]_{odd}$, $[m]_{even}$) the largest integer (resp. the largest odd and the largest even integer) which is less or equal to $m$.\\

\noindent \begin{tabular}[c]{|c|c|l|}
\hline \hline 
\multicolumn{3}{|c|}{\raisebox{-2pt}[5pt][7pt]{\Homc}} \\ 
\hline \hline 
\raisebox{-2pt}[5pt][7pt]{$n=4$} & \raisebox{-2pt}[5pt][7pt]{$m=2k$} & \raisebox{-2pt}[5pt][7pt]{homotopy equivalent to $S^0$} \\
\hline 
\raisebox{-2pt}[5pt][0pt]{$n=2l$} & \raisebox{-2pt}[5pt][0pt]{$m=2sl$} & \raisebox{-2pt}[5pt][0pt]{$2n$ points and $4s-2$ connected components which}\\
\raisebox{-2pt}[5pt][0pt]{$l\neq 2$} & & are homotopy equivalent to $S^1$:\\ 
 & &\raisebox{-2pt}[0pt][10pt]{$\Delta^{1,2}_l,\Delta^{1,2}_{2l},\dots,\Delta^{1,2}_{(2s-1)l}$}\\
\cline{2-3}
 & \raisebox{-2pt}[5pt][7pt]{$m=2k$} & \raisebox{-2pt}[5pt][7pt]{$2(2[\frac{k}{l}]+1)$ connected components, all homotopy}\\
 & {$l\nmid k$} & \raisebox{-2pt}[0pt][10pt]{equivalent to $S^1$: $\Delta^{1,2}_k,\Delta^{1,2}_{k\pm l},\dots,\Delta^{1,2}_{k\pm [k/l]l}$} \\\hline
\raisebox{-2pt}[5pt][7pt]{$n=2l$} & \raisebox{-2pt}{$m=2k+1$} & \raisebox{-2pt}{$\emptyset$} \\
\hline
\raisebox{-2pt}[5pt][0pt]{$n=2l+1$} & \raisebox{-2pt}{$m=sn$} & \raisebox{-2pt}{$2n$ points and $s-1$ connected components $\simeq S^1$:}\\
 & &  \raisebox{-2pt}[0pt][10pt]{$\Delta_n,\Delta_{2n},\dots,\Delta_{(s-1)n}$} \\
\cline{2-3}
 & \raisebox{-2pt}[5pt][0pt]{$m=2k+1$} & \raisebox{-2pt}[5pt][0pt]{$[\frac{m}{n}]_{odd}+1$ components homotopy equivalent to $S^1$:} \\
 & \raisebox{-2pt}[0pt][10pt]{$n \nmid m$} & \raisebox{-2pt}{$\Delta_{(m\pm n[m/n]_{odd})/2 },\Delta_{(m\pm n([m/n]_{odd}-2)/2},\dots,\Delta_{(m\pm n)/2}$ }\\
\cline{2-3} 
 & \raisebox{-2pt}[5pt][0pt]{$m=2k$} & \raisebox{-2pt}[5pt][0pt]{$[\frac{m}{n}]_{even}+1$ connected components, all  $\simeq S^1$:} \\
 & \raisebox{-2pt}[0pt][10pt]{$n\nmid m$} & \raisebox{-2pt}[0pt][10pt]{$\Delta_k,\Delta_{k\pm n},\dots,\Delta_{k\pm \frac{n}{2}[m/n]_{even}}$}  \\
\hline
 \end{tabular}\\
\begin{center}Table 9. \end{center}

\deset
 Since the Euler characteristic of a point, respectively circle, is  equal to 1, respectively 0, we have proven the following claim: 
\begin{cor}\emph{ 
\bdm
\chi(\Hom (C_m,C_n))=\left\{ \begin{array}{ll} 2, & n=4\textnormal{ and } m \textnormal{ is an  even number}; \\ 2n, & \textnormal{when }n \textnormal{ divides }m \textnormal{ and } n\neq 4; \\ 0, & \textnormal{otherwise}. \end{array} \right.
\edm
}
\end{cor} 

On the other hand, we know that if $Y$ is a cell complex, $\chi(Y)=\sum_n(-1)^nc_n$, where $c_n$ denotes the number of $n$-cells of $Y$.\\
Let $n\neq 4$ and let, as in the proof of Theorem \ref{mainth}, $X$ be a connected component of \Homc  with $r$ returning points, where \Homc$\neq \emptyset$. Then dimension of a cell from $X$ belongs to the set $\{0,1,\dots,\min\{r,m-r\}\}$. Let $d\in\{0,1,\dots,\min\{r,m-r\}\}$. We want to find the explicit formula for number $c_d$ of $d$-cells of $X$.  
\begin{itemize}
\item If $r=0$ then $c_0=1$ and $c_d=0$, for $d>0$.
\item Let $r\neq 0$ and $N=n$, if $n$ is odd and $N=\frac{n}{2}$, if $n$ is even. Then $c_0=N\binom{m}{r}$. Let now $0<d\leq\min\{r,m-r\}$.  We define a map $P_d$ which maps any $d$-cell $\eta=(s;s_{i_1},\dots,s_{i_r})$ from $X$ to a $d$-tuple of numbers $(i_{j_1},\dots,i_{j_d})$ such that $j_1<\dots<j_d$ and $s_{i_{j_k}}=i_{j_k}^+$. Also, let $S$ be the set of all $d$-tuples $(\alpha_1,\dots,\alpha_d)$ which satisfy one of the following two conditions:
\begin{eqnarray}
1\leq \alpha_1<\dots<\alpha_d<m \textnormal{ and }\alpha_{j+1}-\alpha_j\geq2\textnormal{, for } j\in[d-1] \label{prva} \\ 
 1< \alpha_1<\dots<\alpha_d=m \textnormal{ and } \alpha_{j+1}-\alpha_j\geq 2 \textnormal{, for } j\in [d-1] \label{druga} 
\end{eqnarray}
It is not hard to see that there exists a bijection between the sets $\{P_d(\eta)\vert \eta\in X,\dim(\eta)=d\}$ and $S$. Since $\vert P_d^{-1}(s)\vert=N\binom{m-2d}{r-d}$, for any $s\in S$, we have that $c_d=N\binom{m-2d}{r-d}\vert S \vert$. But number of $d$-tuples which satisfy (\ref{prva}) is  $\binom{m-d}{d}$ and for (\ref{druga}) is equal to $\binom{m-d-1}{d-1}$. Hence 
\bdm
c_d=N\left[\binom{m-d}{d}+\binom{m-d-1}{d-1}\right]\binom{m-2d}{r-d}=N\frac{m(m-d-1)!}{d!(r-d)!(m-r-d)!}.
\edm
and $\chi(X)=\sum_{d=0}^{\min\{r,m-r\}} (-1)^d N \frac{m(m-d-1)!}{d!(r-d)!(m-r-d)!}$ 
\end{itemize}
Hence we have proven the following formula, for $1\leq r \leq m-1$:
\begin{eqnarray*}
& & \sum_{d=0}^{\min\{r,m-r\}} (-1)^d \frac{(m-d-1)!}{d!(r-d)!(m-r-d)!}= \\ & & \sum_{d=0}^{\min\{r,m-r\}} (-1)^d d \binom{m-d-1}{d-1,r-d,m-r-d}=0.
\end{eqnarray*}

\section{The homotopy type of \textnormal{\Homlk}}

It is easy to see that $\Hom (C_m,L_n)$ is empty if $m$ is an odd integer (see the argument for the fact that \Homc is empty when $n$ is even and $m$ is odd). Hence, from now on we will discuss only the case of \Homlk. Also, since we have already determined the homotopy type of $\Hom (C_{2m},L_3)$ and $\Hom (C_{2m},L_2)$(see the proof of Theorem \ref{l2l3}), we will assume that $n\geq 4$.
\begin{DEF} We say that $i\in [m]$ is \emph{a returning point} of a vertex \temes from \emph{\Homl} if $a_i-a_{[i+1]_m}=1$.
\end{DEF}

 Let $v\in \Hom_0(C_{2m},L_n)$ and let $r$ be the number of its returning points. Then we must have $a_1+2m-2r=a_1$, that is the number of returning points for each vertex must be equal to $m$. \\
 As we did in the previous chapter,  one can prove that there is a bijection between $\Hom_0(C_{2m},L_n)$ and the set $P$, where $P=\{(i;i_1,\dots,i_m)\vert i\in [n],1\leq i_1<\dots<i_m\leq n\textnormal{ and }1\leq i+j-1-2\big\vert\{q\vert i_q<j\}\big\vert\leq n,\textnormal{ for all }j\in[2m]$. Also, we will use the same notations for the cells of \Homl as we did for complexes \Homck.
\begin{rem} Let \temels$=$\temel be a vertex from \Homlk. If for some $l$, $[i_l+1]_{2m}\neq i_{[l+1]_m}$, then $(i;i_1,\dots,i_l+1,\dots,i_m)\in \Hom _0(C_{2m},L_n)$ (we have only replaced $i_l$ with $i_l+1$) if and only if $a_{i_l}=i+i_l+1-2l\neq n$.
\end{rem} 
\begin{lem} \label{conn_compl} Two vertices $u=$\temel and $v=(j;j_1,\dots,j_m)$  of a complex \emph{\Homl} are in the same connected component if and only if $i$ and $j$ have the same parity.
\begin{proof} It is easy to see that, if $i \neq j \mod 2$, those vertices cannot be in the same connected component.\\
Suppose now that $i$ and $j$ have the same parity. Without any loss of generality we can assume that $i\leq j$ and that $j=i+2q$, for some non-negative integer $q$. Like in the proof of Lemma \ref{conn_comp2}, we define equivalence relation $\sim$: for two vertices $x$ and $y$, $x\sim y$ if and only if they lie in the same connected component. \\
Let now, for all  $l\in[n]$, $t^{(l)}=\min \{m,n-l\}$ and let 
\bdm
t^{(l)}_{m-k}=\left\{ \begin{array}{ll} 2m-k, & t^{(l)}\geq 1  \textnormal{ and } k\in\{0,1,\dots,t^{(l)}-1\};\\ 2m-2k+t^{(l)}-1, & m\geq t^{(l)}+1 \textnormal{ and } k\in \{t^{(l)},\dots,m-1\}.  \end{array} \right.
\edm 
If we have in mind the previous remark, it not hard to see the following:
\bdm
\begin{split}
(i;i_1,\dots,i_{r-1},i_m) &\sim (i;i_1,\dots,i_{m-1},t^{(i)}_m)\sim (i;i_1,\dots,t^{(i)}_{m-1},t^{(i)}_m)\sim \cdots \sim \\ &\sim(i;t^{(i)}_1,\dots,t^{(i)}_{m-1},t^{(i)}_m) 
\end{split}
\edm 
It is now clear that, if $i=j$ then $u\sim v$. Suppose that $i<j$. Then, for all $i\leq l<j\leq n$ we have that $t^{(l)}_m=2m$ and $t^{(l)}_1\neq 1$. Hence:
\bdm
\begin{split}
(i;i_1,\dots,i_m) &\sim (i;t^{(i)}_1,\dots,t^{(i)}_{m-1},t^{(i)}_m) \sim (i+2;1,t^{(i)}_1,\dots,t^{(i)}_{m-1})\\ &\sim (i+2;t^{(i+2)}_1,\dots,t^{(i+2)}_{m-1},t^{(i+2)}_m) \sim \dots \sim \\ &\sim (i+2q\,;t^{(i+2q)}_1,\dots,t^{(i+2q)}_m) \sim  (j;j_1,\dots,j_m).
\end{split}
\edm
\end{proof}
\end{lem}
\begin{thm} \emph{\Homl} is homotopy equivalent to two points.
\begin{proof} From Lemma \ref{conn_comp2} we know that \Homl has two connected components, namely $\{(s;s_{i_1},\dots,s_{i_m})\in \Hom (C_{2m},L_n)\vert s \textnormal{ is odd}\}$ and $\{(s;s_{i_1},\dots,s_{i_m})\in \Hom (C_{2m},L_n)\vert s \textnormal{ is even}\}$. Let $X$ be any of these components. We will now prove that $X$ is contractible.  \\
For all $i\in [n]$, define complexes $X_i$ and $\xt$ in the same way as in the proof of Theorem \ref{mainth}. We see that $\xt \subseteq X_i$, $X=\bigcup_{i=1}^n X_i$ and $X_2=X_4=\cdots=X_{2[n/2]}=\emptyset$ or $X_1=X_2=\cdots=X_{2[(n+1)/2]-1}=\emptyset$. \\

The proof that $\xt$ is a strong deformation retract of $X_i$, for non-empty $X_i$, is completely the same as in the proof of the Theorem \ref{mainth}.\\
For a cell $\eta=(i;s_{i_1},\dots,s_{i_r})$ let 
\begin{center}
$F(\eta)= \{i_j\vert i_j\neq m,s_{i_j}=i_j,i_j+1\notin \{i_1,\dots,i_r\}\textnormal{ and }i+i_j+1-2j\neq n\}$
\end{center}
 and $F^+(\eta)= \{i_j\vert s_{i_j}=i_j^+\}$.  The maps $R,R^+:\xt \to \mathbb{N}\,_0$, the set $S$ and the map $\nu$ are  also defined analogously to the corresponding objects in the already mentioned proof. By the previous remark, $\nu$ is well defined.
Again, $(S,\nu)$ is a partial matching on $\mathcal{P}(\xt)$. We will now prove that $\sigma=(i;t^{(i)}_1,t^{(i)}_2,\dots,t^{(i)}_m)$, where $t^{(i)}_1,\dots,t^{(i)}_m$ are defined in the proof of previous lemma, is the only critical cell of our matching. Let $(i;s_{i_1},\dots,s_{i_r})=\xi \in \mathcal{P}\setminus S$. Then we have two cases:\begin{enumerate}
\item $\infty=R(\xi)=R^+(\xi)$.\\ 
Let us first prove that $\infty=R(\sigma)=R^+(\sigma)$:
 \begin{itemize} 
\item If $k\in \{t^{(i)},\dots,m-1\}$, then $t^{(i)}=n-i$, $i+t^{(i)}_{m-k}+1-2(m-k)=i+t^{(i)}=n$ and $t^{(i)}_{m-k}\notin F(\sigma)$.
\item  If $ k\in\{0,1,\dots,t^{(i)}-1\}$, then it is easy to see that either $t^{(i)}_{m-k+1}=t^{(i)}_{m-k}+1$ or, for $k=0$, $t_m^{(i)}=2m$ and $t^{(i)}_{m-k}\notin F(\sigma)$.
\end{itemize}
Hence, $F(\sigma)=\emptyset$ and $\infty=R(\sigma)=R^+(\sigma)$.\\
\noindent Since $R^+(\xi)=\infty$, $\xi$ must be a vertex \temelk.
\begin{itemize} 
\item First we will prove that $i_m=t_m^{(i)}$. Since $R(\xi)=\infty$, $i_m=2m$ or $i_m<2m$ and $i+i_m+1-2m=n$. The second case is possible only if $i=n$ and $i_m=2m-1$. In both cases $i_m=t_m^{(i)}$.
\item Suppose now that $i_{m-k}=t^{(i)}_{m-k}$ for some $k\in [m-1]$. Let us prove that $i_{m-k-1}=t^{(i)}_{m-k-1}$. \\
$\circ$ If $i+t^{(i)}_{m-k}+1-2(m-k)\neq n$, then $k\in\{0,1,\dots,t^{(i)}-2\}$ and then we must have $i_{m-k-1}=t^{(i)}_{m-k}-1=t^{(i)}_{m-k-1}$.\\
$\circ$ If $i+t^{(i)}_{m-k}+1-2(m-k)=n$, then $k\in \{t^{(i)}-1,\dots,m-1\}$ and we must have $i_{m-k-1}=t^{(i)}_{m-k}-2=t^{(i)}_{m-k-1}$.
\end{itemize}
Hence we have proven that $\sigma$ is the only cell with property $\infty=R(\xi)=R^+(\xi)$. 
\item $R^+(\xi)<R(\xi)$:\\
Let  $\xi'=(i;b_{i_1},\dots,b_{i_r})$, where 
\bdm
b_{i_j}=\left\{\begin{array}{ll} s_{i_j}, & i_j\neq R^+(\eta); \\ i_j, & \textnormal{otherwise.} \end{array} \right.
\edm 
Obviously, $\xi'\in S$, $\nu(\xi')=\xi$ and $\xi$ is not a critical cell.
\end{enumerate} 
Hence $\mathcal{P}\setminus (S\cup \nu(S))=\{\sigma\}$. \\

This matching is acyclic (see again the proof of \ref{conn_comp2}), and hence by Proposition \ref{morse} there exists a sequence of elementary collapses leading from $\xt$ to $\{\sigma\}$ and hence $\xt$ is contractible. \\

We have proven that all subcomplexes $X_i$ are contractible.\\
Let $\eta=(A_1,A_2,\dots,A_m)$ be a cell such that $\eta \in X_i\cap X_j$, where $i<j$. Then $(A_1\cup \{i\},A_2,\dots,A_m)\in X_i$ and $(A_1\cup \{j\},A_2,\dots,A_m)\in X_j$. But then, for all $k\in A_2\cup A_m$, it must be $k-i=\pm 1,k-j=\pm1$. This is possible only when $i\leq n-2$ and $j=i+2$ and in that case $A_2=A_m=\{i+1\}$.  We conclude that: 
\bdm
X_i\cap X_j=\left\{ \begin{array}{ll} I \times \widetilde{X}_{k+1}^{m-2}, & \textnormal{for } \vert j-i\vert=2, k=\min\{i,j\}\textnormal{ and }X_i\neq \emptyset\,; \\ \emptyset, & \textnormal{otherwise},  \end{array} \right. 
\edm 
where $\widetilde{X}_{k+1}^{m-2}$ is a subcomplex of $\Hom(C_{2m-2},L_n)$, and $I$ is the unit interval. Since $\widetilde{X}_{k+1}^{m-2}$ is contractible, $X_i\cap X_j$ is also contractible. We also see that $X_{i_1}\cap\dots\cap X_{i_t}=\emptyset$ for $t\geq 3$ (using the similar argument as in case $X_i\cap X_j$).\\
 The family of subcomplexes $\{X_i\}_{i=1}^n$ satisfies the conditions of \cite[Corollary 4G.3, Exercise 4G.4]{Hat} and hence $X\simeq \mathcal{N}(X_i)$. But in this case $\mathcal{N}(X_i)=\{i\vert X_i\neq \emptyset\}\cup \{\{i,i+2\}\vert X_i\neq \emptyset\textnormal{ and }i\leq n-2\}$ and, hence, $\mathcal{N}(X_i)\simeq L_k$, where $k$ is the number of non-empty complexes $X_i$, and $L_k$ is viewed as a 1-dimensional simplicial complex. We conclude that $X$ is contractible, and, hence, \Homl is homotopy equivalent to two points.
\end{proof}
\end{thm}

\end{document}